\documentclass[12pt]{article}
\usepackage{amsfonts}
\usepackage{amssymb}
\usepackage{epsfig}
\def\ve{\varepsilon}

\def\half{{\frac{1}{2}}}

\def\d1e{{\partial _1 e}}
\def\xd1e{{X_{\partial _1 e}}}

\def\calb{{\cal B}}

\def\calg{{\cal G}}

\def\Z{{\mathbb{Z}}}

\def\P{{\mathbb{P}}}
\def\R{{\mathbb{R}}}
\def\Z{{\mathbb{Z}}}

\def\proof{{\bf Proof.}~}

\def\half{{\frac{1}{2}}}

\def\qed{{\hfill$\square$\medskip}}

\def\half{{\frac{1}{2}}}

\def\d1e{{\partial _1 e}}
\def\xd1e{{X_{\partial _1 e}}}

\def\calb{{\mathcal{B}}}

\def\calg{{\mathcal{G}}}

\def\prodb1k{{b_1b_2\cdots b_k}}
\def\proda1k{{a_1a_2\cdots a_k}}

\def\e1k{{e_1,\ldots,e_k)}}
\def\e1n{{e_1,\ldots,e_n)}}
\def\e1star{{e_1*\ldots*e_n)}}
\def\v1n{{v_1,\ldots,v_n}}
\def\x0m{{(x_0,x_1,\ldots,x_m)}}

\def\concx0m{{[x_0,x_1][x_1,x_2]\cdots[x_{m-1},x_m]}}
\def\concalpha1n{{\alpha_1\alpha_2\cdots\alpha_n}}
\def\concbeta1n{{\beta_1\beta_2\cdots\beta_n}}
\def\concbeta1m{{\beta_1\beta_2\cdots\beta_m}}

\newtheorem{theoremsec}{Theorem}[section]

\newtheorem{lemmasec}[theoremsec]{Lemma}

\newtheorem{defsec}[theoremsec]{Definition}

\def\gf{geometrically finite~}

\def\gfg{geometrically finite  group~}

\author{\large Roger C. Alperin and Gennady A. Noskov }

\title{Nonvanishing of algebraic entropy for
geometrically finite groups of isometries of Hadamard manifolds}

\date{}
\begin{document}

\maketitle

\begin{abstract}
We prove that any nonelementary \gfg of
isometries of a pinched Hadamard
manifold has nonzero algebraic entropy in the sense of M. Gromov.
 In other words it has uniform exponential growth.
\end{abstract}

\section{Introduction}
 Given a group $\Gamma $ generated by a finite set $S$ we
denote by $B_S(r)$ the the ball of radius $k$ in the
Cayley graph of $\Gamma $ relative to $S$.
  The {\sf exponential growth rate} or
$\omega (\Gamma ,S)=\lim_{k\rightarrow \infty }$ $^k\sqrt{|B _S(k)|}$ is
well defined (by submultiplicativity). We set
\[
\omega (\Gamma )=\inf
\{
\omega (\Gamma ,S)|S
~{\rm ~is ~a ~finite ~generating~ set~ for~ }\Gamma
\}.
\]

The group $\Gamma $ is of {\sf exponential growth} if
$\omega (\Gamma,S)>1$ for some (and hence for every) generating set $S$ and of
{\sf uniform exponential growth} if $\omega (\Gamma)>1$.
(M. Gromov calls $\ln\omega$ the {\sf entropy} of $\Gamma$ -
the reason is that if $\Gamma$ is the fundamental group of a compact
Riemannian manifold of unit diameter, then $\ln\omega(\Gamma,S)$ is a
lower bound for the topological entropy of the geodesic flow of the manifold
\cite{manning}).
 This notion should not be confused with the one given in \cite{avez},
where the notion of entropy is defined for the random walk on the
group associated to the probability measure on the set of generators.

Clearly uniform exponential growth implies exponential growth.
``On first sight, however, there seems to be no reason for the
converse to be true'' \cite{gromovmetric}.
 Recently and rather spectacularly, J. S. Wilson has given an
example of an exponential growth group which is not uniform
\cite{wilson-nonuniform}; it is a group of automorphisms of a
rooted tree.

The most common way to prove uniformity of
exponential growth is to prove the  {\sf UFS-property}.
 We say that $\Gamma $ satisfies
{\sf UFS-property (=uniformly contains free nonabelian semigroup)}
if there is a constant $n_\Gamma \geq 1$ such that for every
generating set $A$ of $\Gamma $ there exist two elements in
$\Gamma $ of word length $\leq n_\Gamma ,$ freely generating a
free semigroup of rank 2.  The {\sf UFG-property}, obtained from
the previous one by changing free semigroup to free group is
stronger.

Our main result is the following.

\begin{theoremsec}\label{UF}
 Any nonelementary geometrically finite group of isometries
of a pinched Hadamard manifold satisfies the UFG-property and hence
has nonzero entropy.
\end{theoremsec}

We mention a few ideas about the proof of our main result given
here. If $\Gamma$ is cocompact then it is word hyperbolic and the
theorem follows from Koubi's result \cite{koubi} (which was known
before due to M. Gromov \cite{gromovthriller} and T. Delzant
\cite{delzantsous} in the torsionfree case.)  In the main case of
a noncocompact lattice we make use the geometry of {\sf neutered
space} associated to the group.

For a survey on this subject see the book of Gromov and the
articles \cite{grigorsurvey},\cite{harpe-uniform}.
 We mention now some related results.
A group is called {\sf large} if it has a subgroup of finite index
which has a homomorphism onto a free non-abelian group. For
example, in \cite{grunnoskov}, largeness is proven for the
noncocompact lattices in 3-dimensional hyperbolic space.

Also in our article \cite{alperin-noskov-uniform},  partial results were
obtained for subgroups of $\mathrm{SL}_2$ in characteristic zero,
but more complete results  for non-zero characteristic.

\begin{theoremsec}
If a finitely generated subgroup $\Gamma $of $\mathrm{GL}_2(K),K$
a field of nonzero characteristic has exponential growth then it
satisfies the UFG-property and consequently has uniform
exponential growth.
\end{theoremsec}

A. Eskin, S. Mozes and H. Oh have recently given the 'uniform'
proof of Tits' theorem in their beautiful result.

\begin{theoremsec}\cite{EMO}
A finitely generated  linear group over a field of characteristic
zero  has uniform exponential growth if and only if it has
exponential growth.
\end{theoremsec}

\section{Miscellany about $\delta$-hyperbolic spaces}

Let $X$ be a metric space. We denote by $|x-y|$ the distance
between the points $x$ and $y$ in $X$.
 If $X$ is equipped with a basepoint $x_0$ then we shall use the
notation $|x|=|x-x_0|, x\in X.$
 One of definitions of a hyperbolic space  $X$ is given by means of the Gromov
product relative to a base point $x_0$
\begin{equation}
(x\cdot y)=(x \cdot y)_{x_0} = \half
(|x| + |y| - |x-y|), x,y\in X
\end{equation}
as follows: A metric space $X$ is called
{\sf hyperbolic}  \cite{gromovthriller} if there exists a
constant  $\delta \geq 0$ such that for every triple
$x,y,z \in X$ and for every choice of a basepoint the following
holds
\begin{equation}
(x \cdot y) \ge \min ((x \cdot z), (y \cdot z)) - \delta .
\end{equation}

 A {\sf tripod} $T$ is a union of three segments in $\R^2$, which have only
the  origin in common.
 Every geodesic triangle $\Delta=xyz$ in a geodesic metric space $X$
can be mapped onto a tripod so that  the restriction of the map to each side
of $\Delta$ is an isometry.
 We will call a map with these properties a {\sf tripod map}.
 The tripod map always exists and is essentially unique.
 The three points $x',y',z'$ on the sides opposite to $x,y,z$ respectively,
whose image by the tripod map is equal to 0 describe an {\sf
inscribed triangle} and are called the {\sf internal vertices} of
$\Delta$.
 A geodesic triangle $\Delta$ is called $\delta$-{\sf thin} for $\delta\ge 0$,
if the fibers of the tripod map $f:\Delta\rightarrow T$ are of diameter
at most $\delta$.

\begin{lemmasec}\label{thin} Any geodesic triangle in a
$\delta$-hyperbolic space $H$ is $4\delta$-thin.
\end{lemmasec}

\proof See \cite{berkeleynotes} or \cite{CDP}.
\qed

 Let $xyz$ be a  geodesic $\delta$-thin triangle with the inscribed
triangle $x'y'z'$.
 Let $\sigma_y(t),\sigma_z(t)$ be  the arc length parameterizations of the
segments $[x,y],[x,z]$ respectively such that at the moment $t=0$ they
are located in $x$.
 We extend the parameterization for all $t\geq 0$ making the paths
stop when they reach the vertices.
 Then by definition of thinness and lemma \ref{thin} the points
$\sigma_y(t),\sigma_z(t)$ are distance at most $4\delta$ apart until the
moment when  they reach the vertices $y',z'$ of the inscribed triangle, that
is

\begin{equation}
  \label{eq:thin}
|\sigma_y(t)-\sigma_z(t)|\leq 4\delta, 0\leq t\leq T=
\half(|x-y|+|x-z|-|y-z|)=(y\cdot z)_x.
\end{equation}

 We say that $\sigma_y(t),\sigma_z(t)$ $4\delta$-{\sf fellow travel} each
other on the segment $[0,T]$.
 More generally we have the following

\begin{defsec}\label{def-l}
The paths $\sigma(t),\tau(t)$ $\ve$-{\sf fellow travel}
each other on the segment $[t_0,t_1]$ if
$|\sigma(t)-\tau(t)|\leq \ve,~ t_0\leq t\leq t_1.$
 The function $\ell(\sigma(t),\tau(t),\ve), \ve>0 $ is
the supremum of the lengths of all time intervals on which
$\sigma(t),\tau(t)$ $\ve$-fellow travel.
\end{defsec}

\begin{lemmasec}{\sf (Fellow traveller property)}\label{ftp}
  Let $xyz$ be a geodesic triangle in a $\delta -$hyperbolic space $X$ and
$|y-z|=c>0.$
 Let $\sigma_y(t),\sigma_z(t)$ be  the arc length parameterizations of the
the segments $[x,y]$ and $[x,z].$ Then

1) $|\sigma_y(t)-\sigma_z(t)|\leq c+8\delta$ for all $t\geq 0.$

2) For  the reverse parameterizations
$\sigma_y^-(t),\sigma_z^-(t)$ which start from $y,z$ respectively,
we have $|\sigma_y^-(t)-\sigma_z^-(t)|\leq c+8\delta, t\geq 0.$

3) Suppose that $\gamma,\gamma'$ are  geodesics with beginning and
ending distance at most $c>0.$
 We have
 $|\gamma(t)-\gamma'(t)|\leq 16\delta +2c$
for all $t\geq 0.$
\end{lemmasec}

\proof
 1) Consider the inscribed triangle for $xyz$.
 If both $\sigma_y(t)$ and $\sigma_z(t)$
do not reach the vertices of inscribed triangle then
$|\sigma_y(t)-\sigma_z(t)|\leq 4\delta.$
 If they do reach then there are points $y',z'$ on the side
$y,z$ distance at most $4\delta$ from $\sigma_y(t),\sigma_z(t)$
respectively hence  $|\sigma_y(t)-\sigma_z(t)|\leq c+ 8\delta. $

\begin{figure}[!ht]
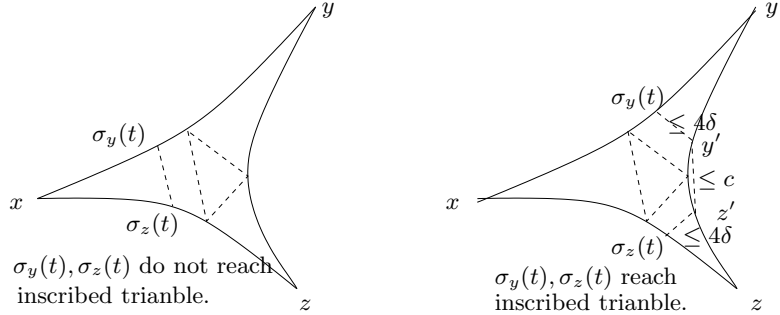
\label{ftpfig}
\centerline{\input ftp.pstex_t}
\caption{Fellow traveller property. Case 1).}
\end{figure}

2) If both $\sigma_y^-(t)$ and $\sigma_z^-(t)$ do not reach the
vertices of inscribed triangle then arguing as in the previous
case we obtain that $|\sigma_y^-(t)-\sigma_z^-(t)|\leq c+8\delta.$
 In contrast  to  the previous case it is possible, that
 $\sigma_y^-(t)$ has reached the vertex of inscribed
triangle but $\sigma_z^-(t)$ has  not.
 If this happens, then
let $z'$ be a point on $[z,x]$ such that
$|\sigma_y^-(t)-x|=|z'-x|$. We have
$
|\sigma_z^-(t)-z'|=
%||z-z'|-t|=
%||z-z'|-|\sigma_y^-(t)-y||=
||z-x|-|y-x|| \leq c,
$
(see Figure 2).
 We conclude that
$|\sigma_y^-(t)-\sigma_z^-(t)|\leq c+4\delta, t\geq 0.$
 Finally  if both $\sigma_y^-(t)$ and $\sigma_z^-(t)$
reach the vertices of inscribed triangle then
again choosing suitable $z'$ (or $y'$) as above
we get  $|\sigma_y'(t)-\sigma_z'(t)|\leq c+ 4\delta.$

\begin{figure}[!ht]
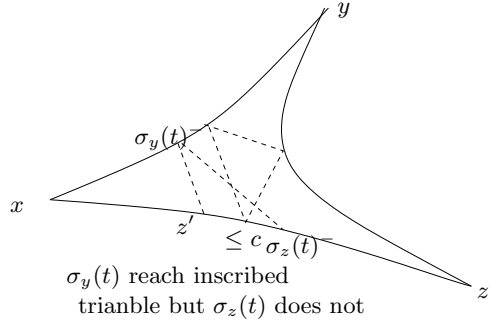
\label{ftp3fig}
\centerline{\input ftp3.pstex_t}
\caption{Fellow traveller property. Case 2)}
\end{figure}

3) Easily  follows applying 1) and 2)
to the triangles $\gamma(0)\gamma(\infty)\gamma'(\infty)$
and $\gamma'(0)\gamma(\infty)\gamma'(\infty)$
respectively.
\qed

 We say that the geodesic triangle $xyz$ has an {\sf obtuse} $y$-{\sf angle}
if either  $y$ is the nearest to  $x$ point on the side $[y,z]$  or
$y$ is the nearest to  $z$ point on the side $[y,x].$

\begin{lemmasec}\label{obtuse}
{\sf (Obtuse-angled triangle lies in $8\delta -$neighbourhood of the side)}
Let $xyz$ be a geodesic triangle with obtuse $y$-angle, then
%$[x,z]\cup [y,z]$ is contained in ${8\delta }$-neighbourhood
%of $[x,y].$
% In particular
$|x-z|\geq |x-y|+|y-z|-16\delta$.
 Furthermore, $(x\cdot z)_y\leq 8\delta.$

\end{lemmasec}

\proof Use the inscribed triangle.
 Say $y$ is nearest to  $x$ on the side $[y,z].$
Since $|z'-x'|\leq 4\delta$ we conclude that $|z'-y|\leq 4\delta$
hence $|y-y'|\leq 8\delta.$
 Thus
$|x-z|= |x-y'|+|y'-z|\geq |x-y|-|y-y'|+|y-z|-|y-y'|
\geq |x-y|+|y-z|- 16\delta$.
 Finally, applying the last inequality we obtain
$(x\cdot z)_y=\half (|x-y|+|z-y|-|z-x|)\leq 8\delta.$

\begin{figure}[!ht]
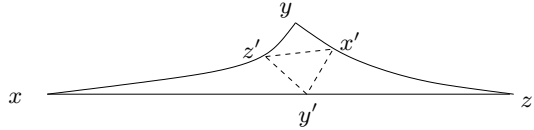
\label{obtusefig}
\centerline{\input obtuse.pstex_t}
\caption{Obtuse triangle inequality.}
\end{figure}

\qed

The next lemmas help us figure out the  structure of geodesic
polygons; they  are proven in greater generality in
\cite{olshanskii-periodic}, see Lemmas 21 and 22.

\begin{lemmasec}\label{remote-gon}
 1) Suppose that  a geodesic $4$-gon
$[x_{1}, x_{2},x_{3},x_{4}]$  satisfies the
following conditions:
$|x_{i+1}-x_{i}| > 180\delta$ for $i =1, 2,3$ and
$(x_{1} \cdot x_{3})_{x_{2}},(x_{2} \cdot x_{4})_{x_{3}}
\leq 14\delta$ for.
 Then the polygonal line $p =x_{1}x_{2}x_{3}x_{4}$ is contained in
the closed $28\delta$-neighborhood of the side
$[x_{1},x_{4}]$.
 In particular,
$|x_{1}-x_{4}| > \sum_1^3 |x_{i+1}-x_i|-168\delta$.

 2) Suppose that  a geodesic $5$-gon $[x_{1}, x_{2},x_{3},x_{4}, x_{5}]$
satisfies the following conditions:
$|x_{i+1}-x_{i}| > 180\delta$ for $i = 1,2,3,4$ and
$(x_{i+2} \cdot x_{i})_{x_{i+1}}\leq 14\delta$ for $i = 1,2,3$.
Then the polygonal line $p =x_{1}x_{2}x_{3}x_{4}x_{5}$ is contained in
the closed $28\delta$-neighborhood of the side
$[x_{5},x_{1}]$.
 In particular,
$|x_{1}-x_{5}| > \sum_1^4 |x_{i+1}-x_i|-168\delta$.
\end{lemmasec}

We have the following variant of the above lemma  for "small" sides
$[x_2,x_3],[x_4,x_5].$

\begin{lemmasec}\label{nearby-gon}
 Suppose that  a
geodesic $5$-gon $[x_{1},x_2,x_3,x_4, x_5]$  satisfies the
conditions
$|x_2-x_3|,|x_4-x_5|\leq 180\delta.$
 Then
$$
|x_{1}-x_{5}|
>
|x_1-x_2|+|x_3-x_4|-360\delta-2\ell_0
,
$$
 where $\ell_0= \ell(\gamma,\gamma',380\delta)$ - the time  of
$380\delta$-fellow travel
of geodesics $\gamma=[x_2,x_1],\gamma'=[x_3,x_4].$

\end{lemmasec}

\proof
First note that $\ell_0\geq (x_1\cdot x_5)_{x_2}.$
  Indeed, the geodesics along the segments $[x_2,x_1],[x_2,x_5]$
~ $4\delta-$fellow travel each other until the moment $t=(x_1\cdot
x_5)_{x_2}$ by the thinness property, see inequality (\ref{eq:thin}).
 Geodesics $[x_2,x_5],[x_3,x_4]$ ~ $376\delta$-fellow travel each
other all the time by the assertion 3) of Lemma \ref{ftp}.
 Thus the segments
$[x_2,x_1],[x_3,x_4]$ ~ $380\delta$-fellow travel each other
 until the moment $t=(x_1\cdot x_5)_{x_2},$ hence
$\ell_0\geq (x_1\cdot x_5)_{x_2}.$
 Finally,
$|x_{1}-x_{5}| =|x_1-x_2|+|x_2-x_5|-2(x_1 \cdot x_5)_{x_2}
\geq|x_1-x_2|+|x_3-x_4|-360 \delta-2(x_1 \cdot x_5)_{x_2}
>
|x_1-x_2|+|x_3-x_4|-360 \delta-2\ell_0
.
$
\qed

\begin{figure}[!ht]
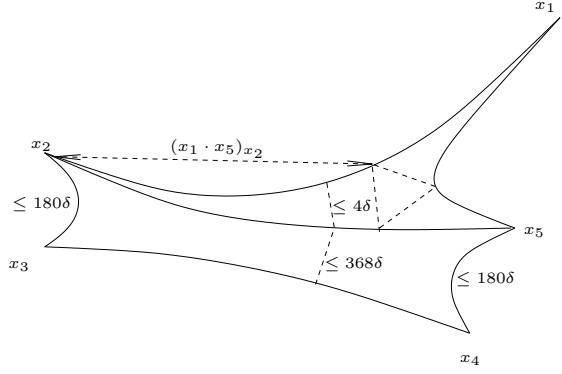
\label{nearby-gon-fig}
\centerline{\input nearby-gon.pstex_t}
\caption{Lemma 2.6: Divergence for 5-gons with two small sides.}
\end{figure}

\section{Uniformly hyperbolic isometries in
groups acting on $\delta-$hyperbolic spaces}

 An isometry $g$ of a $\delta -$hyperbolic space
$X$ is {\sf hyperbolic} if there is a point $x\in X$ such that the
map $n\mapsto g^nx$ from $\Z$ to $X$ is a quasi-isometric
embedding, \cite{CDP}.
This means that $|g^n(x)-x|\geq cn$ for a suitable
positive constant $c$ and every natural $n$.
 The {\sf axis} of $g$ is a  bi-infinite geodesic on which
$g$ acts as a translation.
 In general, a hyperbolic isometry $g$ does not need to possess
an axis, however it does in case  $X$ is a pinched Hadamard
manifold (they are complete CAT(0)-spaces).
 Or $X$ is a Cayley graph of a word
hyperbolic group $\Gamma$ and $g=h^n$ for a constant $n=n(\Gamma)$
and an element $h$ is of infinite order \cite{delzantsous}.

\begin{lemmasec}\label{crithyp} {\rm (\cite{CDP}, lemma 9.2.2.)}
Suppose that for  an isometry $g$ of a
$\delta -$hyperbolic space
$X$ there is a point $x\in X$ such that
\[
|g^2x-x|>|gx-x|+2\delta.
\]
Then $g$ is hyperbolic.
\end{lemmasec}

\begin{lemmasec}\label{hh=h}{\rm (\cite{CDP}, lemma 9.2.3.)}
Suppose that $\delta >0.$ Let $g,h\in Isom X$ be such that for
some $x\in X$ the following holds
\[
\min \{|gx-x|,|hx-x|\}\geq 2(gx\cdot hx)_x+6\delta .
\]
Then $gh,hg$ are hyperbolic.
\end{lemmasec}

Let $S$ be a finite set of isometries of a $\delta $-hyperbolic
space $X$ .  The {\sf size} of $S$ at $x\in X$ is
\[
|S|_x=\max_{s\in S}|sx-x|.
\]

\medskip

The following is a  variant of  Proposition 3.2 from \cite{koubi}
for  actions of a hyperbolic group.

\begin{lemmasec}
\label{short}{\sf (Long generating system gives rise to
a short hyperbolic isometry)}
Let $X$ be a geodesic $\delta $-hyperbolic space and let
$\Gamma $ be a group of isometries of $X,$ generated by a finite set of
isometries $S$.
 Suppose that $|S|_x>100\delta $ for each $x\in X.$
 Then there is a hyperbolic isometry in $\Gamma $ which is a product of at
most 2 isometries from $S$.
\end{lemmasec}

\proof (\cite{koubi}) Suppose that $|S|_x>100\delta $ for each $x\in X$
and that $S$ consists of nonhyperbolic isometries.
 Let $x_0$ be a point of {\sf almost minimal displacement} for $S$
that is
\[
|S|_{x_0}\leq  \inf_{y\in X}|S|_y+ \delta.
\]
 Let $S_0=\{a\in S:|ax_0-x_0|\geq
50\delta \}$ be the {\sf long part} of $S.$
 Fix $a_1\in S$ such that $|a_1x_0-x_0|=\max_{a\in S}|ax-x|>100\delta.$
%Finally put $l=\min \{(a_1x_0\cdot ax_0)_{x_0},a\in S_0\}.$

Case 1): {\it There is $a\in S_0$ such that
$(a_1x_0\cdot ax_0)_{x_0}\leq 20\delta.$}
 Applying the lemma \ref{hh=h} we obtain that $aa_1$ is hyperbolic.

Case 2): {\it $(a_1x_0\cdot ax_0)_{x_0}\geq 20\delta$ for all $a\in S_0$.}
 We will show that this assumption leads to a contradiction of the
almost minimality of $x_0.$
 Namely, we shall show that for the point $x_1$ on the geodesic
$[x_0,a_1x_0]$ such that $|x_1-x_0|=20\delta,$  the following
holds
\[
|S|_{x_1}<|S|_{x_0}-10\delta;
\]
and this contradicts the almost minimality of $x_0.$

 Choose $x_a\in [x_0,ax_0]$ so that $|x_a-x_0|=20\delta,$ then,
since $(a_1x_0\cdot ax_0)_{x_0}\geq 20\delta ,$ by
applying the inequality (\ref{eq:thin}) to
the triple of points $(x_0,a_1x_0,ax_0),$ we have that the
$|x_1-x_a|\leq 4\delta.$
 Consider $x_1$ as a new base point and let $a\in S.$

\begin{figure}[!ht]
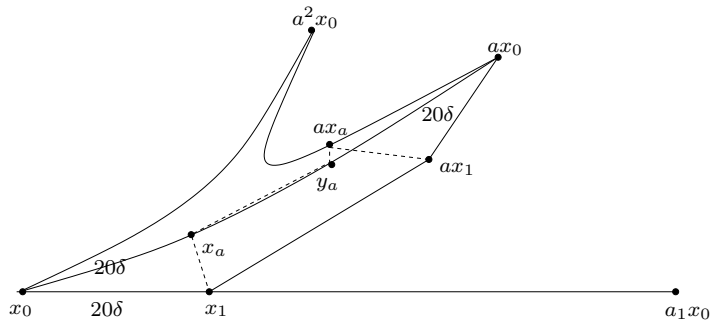
\label{mainlemma}
\centerline{\input main.pstex_t}
\caption{Long generating system gives rise to a uniform hyperbolic isometry.}
\end{figure}

We assert the following:

\begin{itemize}
\item  if $a\in S_0,$ then $|ax_1-x_1|\leq |ax_0-x_0|-28\delta ;$

\item  if $a\notin S_0,$ then $|ax_0-x_0|\leq 90\delta .$
\end{itemize}

It follows from this assertion that
\[
\max_{a\in S}|ax_1-x_1|\leq \max \{\max_{a\in
S}|ax_0-x_0|-28\delta ,90\delta \}<\max_{a\in
S}|ax_0-x_0|-10\delta,
\]
 contradicting the almost minimality of $x_0.$

To prove the assertion,
 consider as a first case when $a\in S_0.$
Choose $y_a\in [x_0,ax_0]$ so that  $|y_a-ax_0|=20\delta .$
We have

\begin{equation}
(a^2x_0\cdot x_0)_{ax_0}=
\half(|ax_0-x_0|+ |a^2x_0-ax_0|-|a^2x_0-x_0|)
\geq \frac 12|ax_0-x_0|-\delta \geq 24\delta .
\end{equation}
 Applying inequality (\ref{eq:thin})
to  the triangle $\Delta (x_0,ax_0,a^2x_0)$ we obtain
$|ax_a-y_a|\leq 4\delta .$
 By assumption $(a_1x_0\cdot ax_0)_{x_0}\geq 20\delta$ hence
$|ax_a-ax_1|=|x_a-x_1|\leq 4\delta .$

We conclude that
\begin{eqnarray*}
|ax_1-x_1| &\leq &|ax_1-ax_a|+|ax_a-y_a|+|y_a-x_a|+|x_a-x_1|\leq \\
&&4\delta +4\delta +(|ax_0-x_0|-40\delta )+4\delta \\
&\leq &|ax_0-x_0|-28\delta .
\end{eqnarray*}

 For the second case when $a\notin S_0$ we have
\begin{eqnarray*}
|ax_1-x_1| &\leq &|ax_1-ax_0|+|ax_0-x_0|+|x_0-x_1|\leq \\
\ 20\delta +50\delta +20\delta &\leq &90\delta
\end{eqnarray*}
\qed

% \begin{lemmasec} {\rm \cite{koubi}, Pr. 3. 2.}
% Let $\Gamma $ be a group acting properly and cocompactly on a
% $\delta$-hyperbolic space $X$.
%  Suppose that $\partial X$ has at least 3 points.
%  Then  there is a constant $n=n(\Gamma,X)$
% depending only on $\Gamma $ and $X$ such that for every finite generating
% set $S$ of $\Gamma $ there is a hyperbolic element $\gamma \in \Gamma $ of
% $S-$length  $\leq n$.
% \end{lemmasec}

% \proof  Let $\calg_0$ be the set of all finite generating sets $S$
% for $\Gamma $ such that $|S|_x>100\delta $ for each $x\in X.$
%  By lemma \ref{short} for any $S$ from $\calg_0$ there is a
% hyperbolic element  in $\Gamma $ which is a product of at most
% 2 isometries from $S$.
%  The remaining generating systems constitute the set
% $\mathcal{G}_1$  of all finite
% generating sets $A$ for $\Gamma $ such that $|X|_x\leq 100\delta $ for some
% $x\in X_\Gamma .$
%  It is easy to see from the assumptions that this set is finite
% up to conjugacy.
%  Since $\Gamma$ is not elementary it contains a hyperbolic
% element.
%  We can write  it (or its conjugate) in generators
% $S\in \mathcal{G}_1$ as a word of uniform length.
% \qed

\section{
Geometrically finite groups of isometries of Hadamard manifolds}

An {\sf Hadamard manifold} is a simply  connected,
complete Riemannian manifold without boundary and with nonpositive sectional
curvatures.
We will assume moreover that $X$ is {\sf pinched}, i.e. all the
sectional curvatures lie between two negative constants:
$-\kappa^2\leq K(X)\leq -1,$
\cite{bowdfinitenessvar}.
 A convenient reference for such manifolds is \cite{BGS}.

A basic fact about Hadamard manifold is that the exponential map
based at any point is injective.
 Thus, any such manifold $X$ is diffeomorphic to $\R^n$.
$X$ can be naturally compactified by adjoining an ideal sphere
$\partial X$ to $X$, so that $\overline{X}=X\cup\partial X$
is homeomorphic to a closed $n$-dimensional ball.
 A pinched Hadamard manifold is a {\sf visibility manifold},
i.e. any two points $x,y\in \overline{X}$ can be joined by a unique
geodesic, which we denote by $[x,y]$.
 When we speak of geodesic as paths, it will be assumed they are parameterized by arc-length.
 Denote by  $|x-y|$  a Riemannian path distance between $x,y\in X.$
 The point $\xi \in \partial X$ can be represented by
a geodesic ray $c_\xi:[0,\infty)\rightarrow X.$ The  function
$h_c(x)=\lim_{t\rightarrow\infty}(|x-c(t)|-t)$ is called the {\sf
horofunction} about $\xi$ associated to the geodesic ray $c$. It
turns out that $h$ is $C^2$ \cite{heintze} and the
norm of its gradient is everywhere equal to 1.
 The level sets of $h_c$ are called
{\sf horospheres} about $\xi.$
 The horospheres form a codimension 1 foliation of $X$
orthogonal to the foliation by bi-infinite
geodesics having one endpoint at $\xi.$
 A set of the form $h_\xi^{-1}[r,\infty)$ for $r\in \R$
is called a {\sf horoball} about $\xi$.
 Such a horoball may alternatively be described as the closure
of the set
$\cup\{\{B_t(\beta(t)):t\in(0,\infty)\}$ where $\beta$ is a geodesic
ray tending to $\xi$ with $\beta(0)\in h_\xi^{-1}(r).$
 In particular horoballs are convex.

 The well known classification of isometries of a hyperbolic
space holds for $X$ too.
  Any isometry $g$ of $X$ extends to a homeomorphism of
$\overline{X}$.
 We shall write $\overline{X}^g$ for the set of fixed points of $g$
in $\overline{X}$.
 Any nonidentical isometry $g$ of $X$ is precisely on of the
following types:

 1) $g$ is {\sf elliptic}, that is $X^g$ is nonempty,

2) $g$ is {\sf parabolic}, that is $\overline{X}^g$ consists of
a single point $p\in \partial X$ and $g$ preserves setwise each horosphere
about $p$,

3) $g$ is {\sf hyperbolic}, that is $\overline{X}^g=\{p,q\},$ where
$p$ and $q$ are distinct points of $\partial X$.
 In this case $g$ translates along the {\sf axis} - the
geodesic connecting $p$ and $q$ and we denote
by $||g||$ the amplitude of this translation
(= {\sf translation length}.

Note that this definition of hyperbolicity is consistent with the
one given before for isometries of a $\delta$- hyperbolic space.

By definition, the group $\Gamma\leq\mathrm{Isom}~X$ is {\sf
elementary} if $\overline{X}^\Gamma\not=\emptyset$ or else
$\Gamma$ preserves bi-infinite geodesic in $\overline{X}.$

\begin{defsec}
\label{gf}{\rm
Let $X$ be a pinched Hadamard manifold and $\Gamma\leq\mathrm{Isom}~X$ is
a subgroup of isometry group, acting properly.
 $\Gamma$ is {\sf geometrically finite} if the  following properties hold:
\begin{itemize}
%\item the $\Gamma$-action is properly discontinuous, i.e.
%      the  $\{g\ |\ g C\cap C\not= \emptyset\}$ is finite for each compact
%      $C\subseteq X$;
\item there is a $\Gamma$-invariant family (possibly empty)
       $\calb=\{B_p\ |\ p\in P\}$ of pairwise disjoint
      closed horoballs in $X$ with $\Gamma\backslash\calb$ finite;
\item there is a non-empty closed convex $\Gamma$-invariant subset
$Z\subseteq X$
such that the  $\Gamma$-action on the {\sf neutered space}
$X_\Gamma=Z- \bigcup_{p\in P}\mathrm{Int}(B_p)$ is cocompact.
\end{itemize}
}
\end{defsec}

The more restrictive class of groups, consisting of  {\sf
lattices},
 corresponds to the case of $Z=X$ in the above definition.

We now return to the discussion of the geometry of $X_\Gamma$.
Each maximal parabolic subgroup $P$ of $\Gamma$ fixes a point at infinity
of $X$ and hence fixes any horoball in $X$ centered at this
point.
 We consider $X_\Gamma = Z - \bigcup_{P\in \P} Int(B_P)$
with the path metric, that is the metric given by lengths of paths,
computed using the standard hyperbolic Riemannian metric.
 $X_\Gamma$ is
complete and locally compact, so it is a geodesic metric space. $\Gamma$
acts on $X$ by isometries, with finite stabilizers, and with compact
quotient.

We call the boundary piece $S_P=X\cap\partial B_P$ that results
from removing $int(B_P)$ a {\sf horosphere of} $X$. Two horospheres
$S_P$ and $S_{P'}$ are have the same image in $M$ if and only if $P$
and $P'$ are conjugate in $\Gamma$.

Since $Z$ is convex, the metric on $Z$ is the restriction of
the metric on $X$.
 In particular, it is a geodesic metric
space, geodesics are unique, and they vary continuously with
choice of endpoint.  For points $x,y\in X$ we will denote
$d_X(x,y)$ their distance apart in $X$ and $d_{X_\Gamma}(x,y)$
their path distance in $X_\Gamma$.
 The path metric on $X_\Gamma$ is at most exponentially distorted with
respect to $d_X;$ namely
$$
d_X(x,y)\leq d_{X_\Gamma}(x,y)\leq\sinh(\kappa^2\cdot d_X(x,y))
$$
for all $x,y\in X$.
 This is proven in \cite{farbcombing} for the case of lattices
but the proof goes equally well for geometrically finite groups.
 Thus, there is a monotonically increasing
positive function $\phi,$ does not depending on $\Gamma$ such that
$\phi (x)\rightarrow \infty$ with $x\rightarrow\infty$ and such
that $d_X>\phi (d_{X_\Gamma}).$

By a result of Bowditch all $P$ are virtually nilpotent
\cite{bowdfinitenessvar} (this follows from the Margulis Lemma).

\section{Uniform hyperbolic elements in groups acting on Hadamard manifolds}

 In the case of noncocompact discrete actions we know of no results
about uniform hyperbolic elements as proven in the previous
section. But we have

\begin{lemmasec}\label{unifhypisometrygf}
Let $\Gamma $ be a nonelementary geometrically finite group acting on a
Hadamard manifold $X.$
 There is neutered space $X_\Gamma$
%(obtained by
%removing from the convex hull of the limit set
%$\Lambda \Gamma$ the disjoint $\Gamma -$
%invariant union of horoballs about parabolic points)
and a  constant $\Delta >0$ such that if $\Gamma $ is generated by
the set $S$ and $|S|_x^{X_\Gamma}>\Delta $
for each $x\in X,$ then there is an hyperbolic isometry in $\Gamma $ which is a
product of at most 2 isometries from $S$.
(The norm is taken relative to the
path metric on $X_\Gamma.$)
\end{lemmasec}

\proof
 Take any $\Delta>0$ such that $\phi (\Delta )>100\delta,$
where $\phi$ is defined in the previous section.
Choose the horoball system so that that for any two of them
$\phi(d_{X_\Gamma}(B_1,B_2))>100\delta.$
 Now take the generating set $S$ so that $|S|_x^{X_\Gamma}>\Delta $
for each $x\in X_\Gamma.$
 Then $|S|_x^{X}>\phi(\Delta)>100\delta, x\in X_\Gamma. $
 To apply the Lemma \ref{short} above it is enough to ensure that
$|S|_x^{X}>100\delta$ for all $x\in X.$
 We have already done this for $x\in X.$
Now suppose $x$ belongs to some horoball $B$ .
 If $Sx\subset B,$ then $\Gamma $
preserves the horoball $B$ and hence is
elementary - contradiction.
 Thus $sx\notin B$ for some $s\in S.$
 The segment $[x,sx]$ crosses the horospheres $\partial B, s\partial B$
in some points $x_1,x_2$ and we have
$|sx-x|_X
\geq |x_2-x_1|_X
\geq \phi(|x_2-x_1|_{X_\Gamma})
\geq\phi(d_{X_\Gamma}(B,sB))>100\delta.$
 Thus, we can apply the Lemma \ref{short}.
\qed

\begin{theoremsec}\label{thm-uniform-hyp}
A geometrically finite group $\Gamma $, acting on a pinched
Hadamard manifold $X$, uniformly contains hyperbolic elements.
\end{theoremsec}

\proof Let $X$ be a corresponding neutered space, obtained by
removing from the convex hull of the limit set the disjoint
$\Gamma$-invariant union of horoballs about parabolic fixed
points. By Lemma \ref{unifhypisometrygf} there is a $\Delta >0 $
such that if $\Gamma $ is generated by $S$ and
$|S|_x^{X_\Gamma}>\Delta$ for each $x\in {X_\Gamma},$ then there
is a hyperbolic isometry in $\Gamma $ which is a product of at
most 2 isometries from $S$.
 (The norm is taken relative to the path metric on $X_\Gamma.$)
 Denote by  $\mathcal{G}_0$ the  set of all finite generating sets $S$
for $\Gamma $ such that  $|S|_x^{X_\Gamma}\leq \Delta $ for some
$x\in X_\Gamma.$
 Since the action on $X_\Gamma $ is proper and cocompact,
the set $\calg_0$ is finite up to conjugacy.
 Denote by  $\calg_1$ the
set of all finite generating sets $S$ for $\Gamma $ do not belonging to
$\calg_0.$
 We want to get a universal bound for the length of a
shortest hyperbolic element in an arbitrary generating system. Firstly it is
easy to do this for $\mathcal{G}_0.$ Indeed, this length is invariant under
conjugation, and $\mathcal{G}_0$ is finite up to conjugacy. Next, if
$|S|_x^{X_\Gamma}>\Delta $ for each $x\in X,$ then by the proposition above
there is a hyperbolic isometry in $\Gamma $ which is a product of not more
than 2 isometries from $S$.
\qed

\section{Separation of axes}

\begin{defsec}\label{def-axis-separ}
{\rm
By an {\sf axis} of an isometry $g$ of a metric space $X$ we
call in isometric copy of $\R$ inside $X$ on which
$g$ acts by translation.
 We say that the group $\Gamma $ of isometries of a metric space
$X$ satisfies $SA-${\sf property (= Separation of axes)} if
for every $\ve>0$ there is $b(\ve)\geq 0$ such that
for any   $g,h\in \Gamma $ possessing the axes and having
the same translation length
either their axes are asymptotic (=fellow travel each other)
or
the following inequality holds:
$\ell(\sigma _g,\sigma_h,\ve )<b(\ve)a,$
where $\sigma _g,\sigma _h$ are the axes of $g,h$ respectively,
and $a$ is the amplitude of translation of $g,h$
on the axes.
 (See \ref{def-l} for  definition of $\ell$).
}
\end{defsec}

\begin{lemmasec}\label{separation}{\sf (Separation of axes)}
If the group $\Gamma $ acts as a  \gf group of
isometries of a pinched Hadamard manifold $X$,
 then $\Gamma$ satisfies the $SA$-property.
\end{lemmasec}

\proof
 For $r>0$ denote by $b_r$ the maximum of
cardinalities of the sets $S$ of elements of $\Gamma ,$ such that
some orbit $Sx,x\in X_\Gamma$ is contained in the ball of radius $r.$
 Since the action of $\Gamma$ on $X_\Gamma$ is cocompact,
$b_r$ is finite.
 We assert that the
$SA$-property is satisfied with $b(\ve)=b_\ve+1.$

 Suppose the contrary then there is
$\ve >0 $ and hyperbolic isometries  $g,h\in \Gamma $
such that $\ell(\sigma_g,\sigma_h,\ve)> (b_\ve+1) a$,
where $\sigma_g,\sigma_h$ are axes of $g,h$ and $a$ is an amplitude
of translation of $g,h$.
% Let $g,h\in \Gamma $ be the isometries acting by
%translation of amplitude $a $ on the axes $\sigma _g,\sigma _h$
%respectively.
 Denote by  $\sigma _g',\sigma _h'$
the subsegments of $\sigma _g,\sigma _h$ respectively such that
$|\sigma_g'(t)-\sigma _h'(t)|\leq \ve ,0\leq t\leq T,$ and
$T\geq (b_\ve+1) a .$
 Note that $\sigma_g$ cannot stay within any horoball $B\in\calb,$
for time longer than $a$ since otherwise $g$ would preserve $B$
and hence would be parabolic.
 Hence, cutting $T$ by amount $a$  we may assume
that $|\sigma_g'(t)-\sigma _h'(t)|\leq \ve ,0\leq t\leq T,$ $T\geq
b_\ve a $ and $\sigma _h'(0)\in X_\Gamma.$
 Changing if necessarily $g,h$ by their
inverses, we may assume that the action of $g,h$ is coherent with the
natural orientations of  $\sigma _g',\sigma _h'$.
 Then $g^i\sigma_g'(0)=\sigma _g'(ia)\in
\sigma _g',h^i\sigma _h'(0)=\sigma _h'(ia)\in \sigma_h',
i=0,1,\ldots ,b_\ve $
and we conclude that
$|g^i\sigma_g'(0)-h^i\sigma _h'(0)|\leq \ve ,i=0,1,\ldots,b_\ve ,$
hence
$|\sigma _g'(0)-g^{-i}h^i\sigma _h'(0)|\leq \ve ,i=0,1,\ldots ,b_\ve .$
 Thus we have $b_\ve +1$ elements
$g^{-i}h^i\sigma _h'(0),i=0,1, \ldots ,b_\ve,$ inside the ball of radius
$\ve $ about  $\sigma _h'(0).$
 All of these elements lie in $X_\Gamma,$ since by our construction
$\sigma _h'(0)\in X_\Gamma !$
 By definition of $b_\ve,$
at least  two of these elements coincide, say
$g^{-i}h^i\sigma _h'(0)=g^{-j}h^j\sigma _h'(0),0\leq i\not=j\leq
b_\ve .$
 We conclude that $g^{i-j}=h^{i-j}$ for some $1\leq i\not= j\leq b_\ve.$
 It follows that $|\sigma_g(n(i-j))-\sigma_h(n(i-j))|\leq \ve$ for all integral
$n$ and hence the axes fellow travel each other.
 Hence they coincide by uniqueness of geodesics.
\qed

\section{Translation discreteness}
 It is known that
if $\Gamma$ is a group of isometries acting properly  and
cocompactly on a convex (in particular CAT(0)) metric space $M$
then $\Gamma$ is translation discrete in the sense that
translation numbers of its nontorsion isometries are bounded away
from zero, \cite{conner-translation-numbers}, Corollary 3.8.

We generalize this to \gf groups

\begin{theoremsec}\label{translation-discrete}
 Suppose that  the group $\Gamma $ acts as a nonelementary \gf
group of  isometries of a pinched Hadamard manifold $X$.
 Then the action is translation discrete in a sense that
the translation numbers of its
hyperbolic isometries are bounded away from zero.
\end{theoremsec}

\proof Let $\calb$ be the horoball system attached to $\Gamma$
acting on $X$ and $X_\Gamma$ - the corresponding neutered space.
 Decreasing  horoballs from $\calb$ we may assume that
that the distances between two distinct horoballs are bounded away
from zero.  Otherwise, there is a sequence of hyperbolic elements
$g_i$ with translation numbers tending to 0.
 Consider first the case when all the axes $A_i$ avoid $\calb.$
 Fix any $x_i\in A_i\subset X_\Gamma,$ then
$||g_ix_i-x_i||\rightarrow 0$ with $i\rightarrow \infty.$
 Since $\Gamma\backslash X_\Gamma$ is compact there is a
compact subset $D\subset X_\Gamma$ such that $\Gamma D=X_\Gamma.$
For every $i$ there is $h_i\in \Gamma,$ so that $h_ix_i\in D.$
Choosing subsequence we may assume that $h_ix_i\rightarrow x\in D.$
 Then we have $|h_ig_ih_i^{-1}x-x| \leq
|h_ig_ih_i^{-1}x-h_ig_ih_i^{-1}h_ix|+
|h_ig_ih_i^{-1}h_ix-h_ix|+
|h_ix-x|=
2|h_ix-x|+|h_ig_ix-h_ix|=
2|h_ix-x|+|g_ix-x|
\rightarrow 0, i\rightarrow \infty.
$
 This contradicts to the properness of $\Gamma$.

 Now consider the case  when there are infinitely
many axes $A_i$ meeting $\calb$ - we may then
assume that all of them meet $\calb$.
 Since  each $g_i$ translates nontrivially each horoball
it visits,  the translation length
of $g_i$ is at least as large as the distance between distinct
horoballs; thus it is bounded away from zero - a contradiction.
\qed

\section{Free subgroups}

Koubi's criterion for free subgroups is as follows. We generalize
 it to our situation.

\begin{lemmasec}\label{freeness-crit}
{\sf (Freeness criterion)}(\cite{koubi}, Lemma 2.4)
Let $g,h$ be the isometries of a
\newline
$\delta -$hyperbolic space $X.$
 Suppose that for a base point $x_0$ the following holds:
%\begin{eqnarray*}
$$
|g^{\pm 1}x_0-h^{\pm 1}x_0|>\max (|gx_0-x_0|,|hx_0-x_0|)+2\delta ,
$$
$$
|g^2x_0-x_0|>|gx_0-x_0|+2\delta ,
$$
$$
|h^2x_0-x_0|>|hx_0-x_0|+2\delta.
$$
%\end{eqnarray*}
 Then $g,h$ freely generate the free group\textit{\ }$F_2$\textit{.}
\end{lemmasec}

\qed

\begin{lemmasec}\label{lemma-uniform-free}
Let $X$ be a $\delta-$hyperbolic space
and let $\Gamma\leq~{\rm Isom}~X $ be a nonelementary
group acting properly on $X$.
 Suppose that the action $\Gamma$ on $X$ is translation
discrete in the sense of Theorem \ref{translation-discrete}
and satisfies the Separation Axes Property.
 Then there is a constant $m=m(\Gamma,X)$ such that for every
hyperbolic  $g_0\in\Gamma$ having axis and for any generating set
$S$ there exists $s\in S$ such that
$g=g_0^m$ and $h=s^{-1}gs$ freely generate a free group of rank 2.
\end{lemmasec}

\proof Since $\Gamma$ on $X$ is translation discrete there is a
constant $c_0>0$ such the translation length of any hyperbolic
element, possessing an axis, is greater than  $c_0.$
 Further, let $b_0=b(372\delta)$ be the constant given by the SA-property.
 We assert that
$m>\frac{1}{c_0}722\delta+2b_0$ is big enough
to satisfy conclusion of the lemma.
 Suppose  $g_0\in \Gamma$ has an axis $\sigma_g$ and
let $S$ be a finite generating system for $\Gamma$.
 Since $\Gamma $ is nonelementary, there exists $s\in S$ not
fixing the pair $\{\sigma(-\infty),\sigma(\infty)\}.$
 Then the axes $\sigma_{g_0},\sigma_{h_0}$ of
$g_0,h_0=sg_0s^{-1}$ have disjoint limit sets.
 By Separation Axes Property

 \begin{equation}
\ell_0=\ell(\sigma _{g_0},\sigma_{h_0},380\delta )\leq b_0||g_0||.
 \end{equation}

 By the choice of $m$ the translation length $m||g_0||$ of
$g=g_0^m$ is greater than $722\delta+2b_0||g_0||$.
 Reparameterize the geodesics so that at the
moment $t=0$ they are located at the points
$x_0\in\sigma_g=\sigma_{g_0},y_0\in\sigma_h=\sigma_{h_0}$, which
are the closest ones among the points of these axes.

We shall verify conditions of the criterion of freeness
\ref{freeness-crit}. We start with the condition
$$
|gx_0-hx_0|
>
\max
\{|gx_0-x_0|,|hx_0-x_0|\}+2\delta.
$$

{\sf Remote case:}
 Suppose that $|x_0-y_0|=|hx_0-hy_0|>180\delta.$
 The 5-gon $[gx_0,x_0,y_0,hy_0,hx_0]$ satisfies assumptions
of Lemma \ref{remote-gon}.
 Indeed, the lengths of the sides $[gx_0,x_0], [hy_0,hx_0]$
are at least $180\delta$ each.
  Angles at the points $x_0,y_0,hy_0$ are all obtuse,
since, for example, $x_0$ is the nearest to $y_0$
point on  $\sigma_g$.
 Hence by Lemma \ref{obtuse}
 $(gx_0\cdot y_0)_{x_0}<8\delta\leq 14\delta.$

We conclude that

$$
|gx_0-hx_0|
>
|gx_0-x_0|+|hy_0-y_0|+2|x_0-y_0|-168\delta
\geq
$$
$$
|gx_0-x_0|+|hx_0-x_0|-168\delta.
$$

Applying  lemma \ref{remote-gon} to the 4-gon
$[x_0,y_0,hy_0,hx_0]$ we obtain

$$
|hx_0-x_0|>2|x_0-y_0|+|hy_0-y_0|-168\delta > 180\delta,
$$

hence, combining with the above we obtain

$$
|gx_0-hx_0|
>
\max
\{|gx_0-x_0|,|hx_0-x_0|\}+2\delta.
$$

\begin{figure}[!ht]
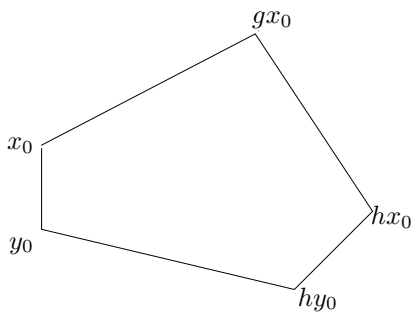
\label{remotefig}
\centerline{\input remote.pstex_t}
\caption{Geodesic 5-gon.}
\end{figure}

{\sf Nearby case:} Suppose that $|x_0-y_0|=|hx_0-hy_0|\leq
180\delta.$
 Then by lemma \ref{nearby-gon} applied to the 5-gon
$[gx_0,x_0,y_0,hy_0,hx_0]$ we have

$$
|gx_0-hx_0|
>
|gx_0-x_0|+|hy_0-y_0|-360\delta -2b_0||g_0||\geq
$$
(since both summands below are
$\geq 722\delta+2b_0||g_0||$)
$$
\max \{|gx_0-x_0|,|hy_0-y_0|\}+362 \delta \geq
$$
(since $|hy_0-y_0|\geq |hx_0-x_0|-360\delta)$)
$$
 \max \{|gx_0-x_0|,|hx_0-x_0|\}+2\delta.
$$

It remains to verify the conditions $
|g^2x_0-x_0|>|gx_0-x_0|+2\delta , |h^2x_0-x_0|>|hx_0-x_0|+2\delta.
$

The first is clear since $|g^2x_0-x_0|=2|gx_0-x_0|.$ For the
second, since $x_0$ does not lie on the axis of $h$, we shall
consider the remote and nearby cases.
 If $|x_0-y_0|=|hx_0-hy_0|>180\delta$
then applying  Lemma \ref{remote-gon} to the 4-gon
$[x_0,y_0,h^2y_0,h^2x_0]$ we obtain

$$
|h^2x_0-x_0|>2|x_0-y_0|+|h^2y_0-y_0|-168\delta
=
2|x_0-y_0|+2|hy_0-y_0|-168\delta
$$

$$
>360\delta+|hy_0-y_0|-168\delta
>
|hy_0-y_0|+2\delta.
$$
 In the nearby case we have $|x_0-y_0|=|hx_0-hy_0|\leq 180\delta$
and then
$$
|h^2x_0-x_0|
\geq
|h^2y_0-y_0|- 2|x_0-y_0|
\geq
2|hy_0-y_0| - 360\delta
>
|hy_0-y_0|+ 360\delta.
$$
\qed

{\bf Proof of the theorem \ref{UF}}.
 Suppose that $\Gamma$ is a nonelementary geometrically finite group
of isometries of a pinched Hadamard manifold $X$.
 Denote by  $X_\Gamma$ the corresponding neutered space.
 Since $\Gamma$ acting on $X_\Gamma$ is proper and cocompact it satisfies
Separation Axes Property. By Lemma \ref{separation} $\Gamma$
acting on $X$ satisfies Separation Axes Property too.
 Moreover by Theorem \ref{thm-uniform-hyp}, $\Gamma$
uniformly contains hyperbolic elements.
 Hence by lemma \ref{lemma-uniform-free} $\Gamma$ uniformly
contains  a free group of rank two. \qed

%\bibliographystyle{alpha}
%\bibliography{aut}

\vspace{15mm}

\newcommand{\etalchar}[1]{$^{#1}$}
\def\cprime{$'$} \def\cprime{$'$}

 \begin{flushleft}
 Roger C. Alperin:\\
 Dept. of Mathematics\\
 San Jose State University\\
 San Jose, CA 95192 USA\\
 E-mail: alperin@math.sjsu.edu
 \end{flushleft}

 \begin{flushleft}
 Gennady A. Noskov:\\
 Institute of Mathematics\\
 Russian Academy of Sciences\\
  Pevtsova 13, Omsk, 644099, Russia\\
 and\\
 Fakult\"{a}t f\"{u}r Mathematik \\
 Universit\"{a}t Bielefeld\\
 Postfach 100131,\\
 D-33501 Bielefeld, Germany\\
 E-mail: noskov@mathematik.uni-bielfeld.de

 \end{flushleft}

 \end{document}